\documentclass{ifacconf}

\usepackage{graphicx}      
\usepackage{natbib}        
\usepackage{scalerel}       

\usepackage{graphics} 
\usepackage{epsfig} 
\usepackage{amsmath} 
\usepackage{amssymb}  

\newcommand{\xpos}{\text{x}}
\newcommand{\xvel}{\text{v}}
\usepackage{xcolor} 
\usepackage{pgfplots} 
\newtheorem{theorem}{Theorem}
\newtheorem{remark}{Remark}
\newtheorem{lemma}[theorem]{Lemma}

\newtheorem{assumption}{Assumption}
\newtheorem{definition}{Definition}
\newtheorem{Problem}{Problem}
\usepackage{tikz}
\usetikzlibrary{shapes.geometric,arrows}
\usepackage[disable]{todonotes}
\usepackage[hyphens]{url}

\pgfplotsset{compat=1.17}
\begin{document}

\begin{frontmatter}

This is an extended version (author's copy) supplementing a paper to be presented at IFAC 2023.
\title{Gradient-based Cooperative Control of quasi-Linear Parameter Varying Vehicles \\ with Noisy Gradients\thanksref{footnoteinfo}} 

\thanks[footnoteinfo]{This work was partially funded by the German Research Foundation (DFG) within
	their priority programme SPP 1914 Cyber-Physical Networking.}

\author[address]{Adwait Datar}
\author[address]{Antonio Mendez Gonzalez}
\author[address]{Herbert Werner}

\address[address]{Institute of Control systems, Hamburg University of Technology,
	Ei{\ss}endorfer Str. 40, 21073 Hamburg, Germany.\\
       {\tt\small \{adwait.datar, antonio.mendez, h.werner\}@tuhh.de}}

\begin{abstract}                
This paper extends recent results on the exponential performance analysis of gradient based cooperative control dynamics using the framework of exponential integral quadratic constraints ($\alpha-$IQCs).
A cooperative source-seeking problem is considered as a specific example where one or more vehicles are embedded in a strongly convex scalar field and are required to converge to a formation located at the minimum of a field.
A subset of the agents are assumed to have the knowledge of the gradient of the field evaluated at their respective locations and the interaction graph is assumed to be uncertain.
As a first contribution, we extend earlier results on linear time invariant (LTI) systems to non-linear systems by using quasi-linear parameter varying (qLPV) representations.
Secondly, we remove the assumption on perfect gradient measurements and consider multiplicative noise in the analysis.
Performance-robustness trade off curves are presented to illustrate the use of presented methods for tuning controller gains.
The results are demonstrated on a non-linear second order vehicle model with a velocity-dependent non-linear damping and a local gain-scheduled tracking controller. 
\end{abstract}


\end{frontmatter}

\section{Introduction}\label{sec:introduction}
Gradient-based forcing terms occur in many well known cooperative control dynamics such as formation control (\cite{Fax.2004}),  flocking (\cite{OlfatiSaber.2006}), extremum-seeking (\cite{michalowsky2016extremum}, \cite{khong2014multi}), etc.
One such application scenario is that of cooperative source-seeking where the main objective is for a group of autonomous vehicles (or a single one), positioned at arbitrary locations within a scalar field (e.g. an oil spill, see \cite{senga2007development}), to locate the minimum/maximum of such field (a.k.a. the source).
A common technique to locate this source is to use the gradient of the field as the indicator of the direction in which the vehicles should move.
Such dynamics have been studied in the literature.
\cite{michalowsky2016extremum} study gradient-based and gradient-free techniques to locate the extremum, where linear matrix inequalities (LMIs) are used to optimize the controller gains.
\cite{Datar.2022} study gradient-based cooperative source-seeking problems for linear time invariant (LTI) and linear parameter varying (LPV) vehicles with exact gradients where a controller is assumed to be designed and the focus is on analysis.
By means of the so called $\alpha$ integral quadratic constraints ($\alpha-$IQCs), LMIs are derived, whose feasibility guarantee exponential convergence of a group of vehicles under formation to the source, with rate $\alpha$ for a range of possible communication networks and scalar fields.
The present paper extends these results by considering non-linear vehicles models and imperfect gradient measurement modeled by deterministic multiplicative noise. 
The non-linear models are represented as quasi-LPV (qLPV) models such that we can deploy established controller synthesis tools (see \cite{Hoffmann.2015} for an extensive survey) to design tracking controllers.
LMIs are derived to guarantee that a group of vehicles locate the source of a scalar field under a predefined formation, with guaranteed exponential rate $\alpha$ for a predefined noise level $\delta$.
Similar analyses have appeared in the literature in the context of analyzing first-order optimization algorithms.
\cite{hu2021analysis} study the stochastic gradient descent with a noise model that includes additive and multiplicative components.
\cite{van2021speed} extend the study to momentum methods in the presence of additive noise.
Although the application studied in the above papers is reasonably different, the analysis we present can be seen as the continuous-time higher-order analogues of these results.

The main contributions of this paper can be summarized as follows:
\begin{enumerate}
    \item Performance analysis results developed by \cite{Datar.2022} which assume perfect gradient measurements are extended by explicitly considering multiplicative noise with a known bound (see Theorem \ref{theom:theorem_perf_analysis_noisy}).
    \item A local result directly applicable to qLPV systems is given which opens the doors for considering non-linear vehicle models (see Theorem \ref{theom:theorem_perf_analysis_noisy_qlpv}).
    \item Minor adjustments to the decomposition result from \cite{Datar.2022} for LTI systems are made to accommodate applications to qLPV systems with heterogeneous scheduling (see Lemma \ref{lemm:decomposition}).
    \item Numerical examples are provided that demonstrate how the developed theoretical results could be used to tune controllers for specific noise levels.
\end{enumerate}




\subsection*{Notation}\label{sec:notation}
The condition number of a matrix $X$ is denoted as $\textnormal{cond}(X)$.
For any $x\in \mathbb{R}^n$, let $\textnormal{diag}(x)$ denote the diagonal matrix formed by placing the entries of $x$ along the diagonal. 
For block matrices, we use $*$ to denote required entries to make the matrix symmetric. 
Let $\mathbf{0}$ and $\mathbf{1}$ denote the vectors or matrices of all zeroes and ones of appropriate sizes, respectively. 
Let $I_d$ be the identity matrix of dimension $d$ and we remove the subscript $d$ when the dimension is clear from context. 
Let $\otimes$ represent the Kronecker product.
Let $\mathcal{S}(m,L)$ denote the set of continuously differentiable functions $f$ which are strongly convex with parameter $m$, and have Lipschitz gradients with parameter $L$ for some given $0<m \leq L$, i.e., 
\begin{equation*}
    m||y_1-y_2||^2 \leq  (\nabla f(y_1)-\nabla f(y_2))^T(y_1-y_2) \leq L ||y_1-y_2||^2
\end{equation*}
holds for all $y_1,y_2$.
Note that for a fixed value of $m$, increasing $L$ enlarges the set $\mathcal{S}(m,L)$.
The set of vector valued functions which are square-integrable over $[0,T]$ for any finite $T$ is denoted by $\mathcal{L}_{2e}[0,\infty)$.
We use 
$\left[\begin{array}{c|c}
\mathcal{A}     &  \mathcal{B}\\
\hline
\mathcal{C}     &  \mathcal{D}
\end{array}\right]$
to represent an LTI system with state-space realization given by matrices $\mathcal{A},\mathcal{B},\mathcal{C}$ and $\mathcal{D}$.
For an ordered set of vectors $(x_1,x_2,\hdots,x_N)$, let the vector formed by stacking these vectors be denoted by $x=\left[x_1^T \; x_2^T \; \hdots \; x_N^T \right]^T$.
Let $\textnormal{blkdiag}(X_1,X_2,\cdots,X_N)$ denote a block-diagonal matrix formed by placing the matrices $X_1$, $X_2$ and so on, as the diagonal blocks.  
Let $\hat{X}$ denote the matrix $I_N \otimes X$ and let $X_{(d)}$ denote the matrix $X \otimes I_d$ for any matrix $X$.
In the context of LPV systems, for an ordered set of parameters $(\rho_1,\rho_2,\hdots,\rho_N)$ and a parameter dependent matrix $X(\rho_i)$, we use $\hat{X}(\rho)$ to denote $\textnormal{blkdiag}(X(\rho_1),X(\rho_2),\cdots,X(\rho_N))$.  
For $N$ identical LPV systems 
$G(\rho_i)=\left[\begin{array}{c|c}
	A(\rho_i)     &  B(\rho_i)\\
	\hline
	C(\rho_i)     &  D(\rho_i)
\end{array}\right]$ parameterized by $\rho_i$ for $i \in \{1,2,\cdots,N\}$, the notation $\hat{G}(\rho)$ represents an LPV system given by
$\left[\begin{array}{c|c}
	\hat{A}(\rho)     &  \hat{B}(\rho)\\
	\hline
	\hat{C}(\rho)     &  \hat{D}(\rho)
\end{array}\right].$ 
\section{Problem Setup}\label{sec:problem_setup}
Consider a source-seeking scenario where $N$ vehicle agents moving in $\mathbb{R}^d$ space (typically $d\in\{1,2,3\}$) are embedded in an underlying differentiable scalar field $\psi: \mathbb{R}^d \xrightarrow[]{}\mathbb{R}$ which satisfies the following assumption.
\begin{assumption}\label{assum:psi_Sml}
Let $\psi \in \mathcal{S}(m_{\psi},L_{\psi})$ and let $y_{\textnormal{opt}}$ minimize $\psi$, i.e., $\psi(y)\geq \psi(y_{\textnormal{opt}}) \, \forall y \in \mathbb{R}^d$ and $\nabla \psi(y_{\textnormal{opt}})=0$.
\end{assumption}
The interconnections between the vehicle agents are modeled with an undirected graph $\mathcal{G}=(\mathcal{V},\mathcal{E})$ such that each vertex $i \in \mathcal{V}$ represents a vehicle and vehicles $i$ and $j$ communicate if and only if $(i,j) \in \mathcal{E}$.
It is assumed that a non-empty subset $\mathcal{V}_l \subseteq \mathcal{V}$ of informed agents (agents with additional information) have access to the local gradient $\nabla \psi$ evaluated at their respective positions.
The dynamics of the these informed agents can therefore be augmented with a forcing term in the direction of the negative gradient that drives them towards the source.
The following assumption on the connectivity of the graph is made, the necessity of which, is discussed later.
\begin{assumption}\label{assm:path_to_leaders}
For every node $i \in \mathcal{V}$, there is a node $j \in \mathcal{V}_l$ such that $\mathcal{G}$ contains a path from $i$ to $j$.
\end{assumption}

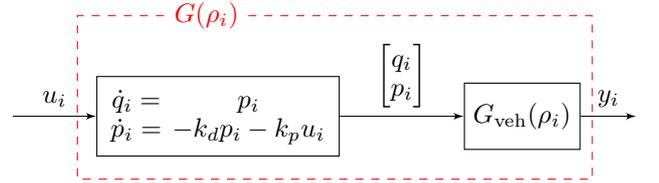
\begin{figure}[h!] 	
 	\centering
 	\tikzstyle{block} = [draw, rectangle, 
    minimum height=3em, minimum width=6em]
\tikzstyle{sum} = [draw,  circle, node distance=1cm]
\tikzstyle{input} = [coordinate]
\tikzstyle{output} = [coordinate]
\tikzstyle{pinstyle} = [pin edge={to-,thin,black}]

%
\scalebox{1}{
\begin{tikzpicture}[auto, node distance=1.5cm,>=latex']
\node [input, name=input] {};
\node [block, right of=input,node distance=2.7cm](Flockcontrol) {$\begin{array}{cc}
   \Dot{q}_i =& p_i \\
   \Dot{p}_i =&-k_d p_i-k_p u_i
\end{array}$};
\node [block, right of=Flockcontrol,node distance=4cm,minimum height=2.5em,minimum width=3em](velocityloop) {$G_{\textnormal{veh}}(\rho_i)$};
\node [output, right of=velocityloop,node distance=1.5cm](output) {};
\node[text=red, above left= 5mm and -2cm of Flockcontrol] (veh) {$G(\rho_i)$};
\draw[red, dashed] (veh.east)-|([xshift=-5.8mm]output.west)|-([yshift=-3mm]Flockcontrol.south)-|([xshift=0.85cm]input.east)|-(veh.west);
\draw [draw,->] (input) -- node {$u_i$} (Flockcontrol);
\draw [draw,->] (Flockcontrol) -- node {$\begin{bmatrix}q_i\\p_i\end{bmatrix}$} (velocityloop);
\draw [draw,->] (velocityloop) -- node[name=outmid] {$y_i$} (output);
\end{tikzpicture}
}
 	\caption{Local control architecture.}
	\label{fig:plant G}	
 \end{figure}

We assume that a local tracking controller has been designed and the closed-loop dynamics of the $i^{\textnormal{th}}$ vehicle agent with desired reference position $q_i(t)$, desired reference velocity $p_i(t)$ can be described by an LPV system $G(\rho_i)$ along with a compact set $\mathcal{P}\subseteq \mathbb{R}^{n_{\rho}}$ such that the function $\rho_i:[0,\infty) \rightarrow \mathcal{P}$ captures the time-dependence of the model parameters.
These dynamics for a given initial condition $x_i(0) \in \mathbb{R}^{n_x}$ can be described by 
\begin{equation}\label{eq:vehicle_dyn}
    \begin{split}
        \Dot{x}_i(t)&=A(\rho_i)x_i(t) + \begin{bmatrix} B_q(\rho_i) & B_p(\rho_i) \end{bmatrix} \begin{bmatrix} q_i(t)\\p_i(t)\end{bmatrix}, \\
        y_i(t)&=C(\rho_i)x_i(t).
    \end{split}
\end{equation}
where the dependence of $\rho_i$ on $t$ is suppressed for compactness, i.e., $\rho_i$ denotes $\rho_i(t)$.
We propose to augment these closed-loop vehicle dynamics shown in Fig. \ref{fig:plant G} by the second order dynamics ,
\begin{equation}\label{eq:vir_vehicle_dyn}
    \begin{split}
        \Dot{q}_i(t)&=p_i(t),\\
        \Dot{p}_i(t)&=-k_d \cdot p_i(t) - k_p \cdot u_i(t),
    \end{split}
\end{equation}
where $u_i(t) \in \mathbb{R}^d$ denotes external input, $q_i(0)=Cx_i(0)$, $p_i(0)=0$ such that $q_i(t)$ and $p_i(t)$ are fed as inputs to dynamics \eqref{eq:vehicle_dyn}.

The overall dynamics with state $\eta_i=\left[x_i^T \; q_i^T\; p_i^T\right]^T$ and suitable initial condition $\eta_i(0)$ can be represented by
\begin{equation} \label{eq:sys_dyn_G}
    \begin{split}
        \Dot{\eta}_i(t)&=A_G(\rho_i)\eta_i(t) + B_G(\rho_i) u_i(t), \\ 
        y_i(t)&=C_G(\rho_i) \eta_i(t), \\
    \end{split}
\end{equation}
where
\begin{equation*}
\left[\begin{array}{c|c}
A_G(\rho_i)     &  B_G(\rho_i)\\
\hline
C_G(\rho_i)     &  D_G(\rho_i)
\end{array}\right]
=
\left[\begin{array}{ccc|c}
   A(\rho_i)            & B_q(\rho_i) & B_p(\rho_i) & \mathbf{0}\\
   \mathbf{0}   & \mathbf{0} & I_d & \mathbf{0}\\
   \mathbf{0}  & \mathbf{0} & -k_d I_d & -k_p I_d\\
   \hline
   C(\rho_i)  & \mathbf{0} & \mathbf{0} & \mathbf{0}\\
\end{array}
\right].
\end{equation*}
We restrict our attention to vehicle models and therefore require integral action in our state-space models. This is captured in the following assumption.
\begin{assumption}\label{assm:equilibrium}
For any $y_{i*} \in \mathbb{R}^d$, there exists an equilibrium $\eta_{i*}$ such that 
\begin{equation} \label{eq:equilibrium}
    \begin{split}
        \mathbf{0}&=A_G(\rho_i(t))\eta_{i*}, \\ 
        y_{i*}&=C_G(\rho_i(t)) \eta_{i*} \\
    \end{split}
\end{equation}
hold for any trajectory $\rho_i:[0,\infty) \rightarrow \mathcal{P}$.
\end{assumption}
\begin{remark}
If the local tracking controller is designed to have zero steady-state error for step position references, i.e., $C(\rho_i(t))A(\rho_i(t))^{-1}B_q(\rho_i(t))=I_d$ for any trajectory $\rho_i:[0,\infty) \rightarrow \mathcal{P}$, Assumption \ref{assm:equilibrium} is satisfied with $y_{i*}=C(\rho_i(t))x_{i*}=q_{i*}$.
\end{remark}
\begin{remark}
For LTI systems, Assumption \ref{assm:equilibrium} reduces to the condition that $(A_G, C_G)$ is detectable (See Theorem 2.1 from \cite{Scherer.12022021}) which is equivalent to $\left(\begin{bmatrix}A &B_q\\\textbf{0}&\textbf{0}\end{bmatrix},[C\,\,\mathbf{0}]\right)$ being detectable.
\end{remark}
With the notation introduced in Section \ref{sec:notation}, let $\eta(t)$, $u(t)$, $y(t)$ and $\rho(t)$ be obtained by stacking the states, inputs, outputs and parameters of the agents, respectively.
The full system can be described by
\begin{equation} \label{eq:sys_dyn_G_big}
    \begin{split}
        \Dot{\eta}(t)&=\hat{A}_G(\rho)\eta(t) + \hat{B}_G(\rho) u(t), \quad \quad \eta(0)=\eta_{0},\\
        y(t)&=\hat{C}_G(\rho) \eta(t), \\
    \end{split}
\end{equation}
where $\rho:[0,\infty) \rightarrow \mathcal{P}^N$ captures the time-dependence of the model parameters for the complete system.
Define $u_{\psi}$ by stacking up
\begin{equation}\label{eq:forcing_term}
    \begin{split}
        u_{\psi_i}(t)&=
        \begin{cases}
        \nabla \psi(y_i(t)), & \text{if}\ i\in \mathcal{V}_l, \\
        0 & \text{otherwise}.
        \end{cases}
    \end{split}
\end{equation}
Standard formation control dynamics with the formation reference $r$ (See \cite{Fax.2004}) and additional forcing term on the informed agents can be represented by
$u=\mathcal{L}_{(d)}(y-r)+u_{\psi}.$
The overall closed-loop system is now described by
\begin{equation} \label{eq:sys_dyn_G_closeloop}
    \begin{split}
        \Dot{\eta}(t)&=\hat{A}_G(\rho)\eta(t) + \hat{B}_G(\rho) u(t), \quad \quad \eta(0)=\eta_{0},\\
        y(t)&=\hat{C}_G(\rho) \eta(t), \\
        u(t)&=\mathcal{L}_{(d)}(y-r)+u_{\psi}.
    \end{split}
\end{equation}
In order to describe the dynamics \eqref{eq:sys_dyn_G_closeloop} as a standard robust control problem as in the work of \cite{Datar.12102021} with the scalar field and communication graph modeled as an uncertainty, define a function $f:\mathbb{R}^{Nd}\rightarrow \mathbb{R}$ as follows.
\begin{definition}\label{defn:f}
For a given graph $\mathcal{G}$ of order $N$ (with its corresponding Laplacian $\mathcal{L}$), the set of informed agents $\mathcal{V}_l$, a scalar field $\psi$ and a given formation reference vector $r \in \mathbb{R}^{Nd}$, define a function $f:\mathbb{R}^{Nd}\rightarrow \mathbb{R}$ by
\begin{equation}\label{eq:defn_f}
    f(y)=\frac{1}{2}(y-r)^T(\mathcal{L} \otimes I_d)(y-r) + \sum_{i \in \mathcal{V}_l} \psi(y_i).
\end{equation}
\end{definition}
Observe that $u(t)=\nabla f (y(t))$.
So the overall dynamics can be described by
\begin{equation} \label{eq:sys_dyn_G_hat}
    \begin{split}
        \Dot{\eta}(t)&=\hat{A}_G(\rho)\eta(t) + \hat{B}_G(\rho) u(t), \quad \quad \eta(0)=\eta_{0},\\
        y(t)&=\hat{C}_G(\rho)\eta(t), \\
        u(t)&=\nabla f (y(t)).
    \end{split}
\end{equation}
Imperfect measurements of inter-agent distance $\mathcal{L}_{(d)}(y-r)$ and field-gradient measurements $\nabla \psi$ can be modeled by a multiplicative noise corrupting the overall measurement of $\nabla f (y(t))$, i.e. the control input $u(t)=\nabla f (y(t))$ exerted on the system can be replaced by 
\begin{equation}
u_n(t)=(I+\Delta)\cdot u(t)=u(t)+\Delta\cdot u(t)=u(t)+e(t).  \label{eq:noise_condition}
\end{equation}
The knowledge about noise is assumed to be available in the form of the following assumption.
\begin{assumption} \label{assum:noise}
For a non-negative constant $\delta$ and arbitrary trajectories $y(t)$ and $e(t)$, let $||e(t)||\leq \delta ||\nabla f(y(t))||$, i.e., 
		\begin{equation}
			\begin{bmatrix}
				e(t)\\\nabla f(y(t))
			\end{bmatrix}^T
			\begin{bmatrix}
				-I & 0\\ 0 & \delta^2I
			\end{bmatrix}
			\begin{bmatrix}
				e(t)\\\nabla f(y(t))
			\end{bmatrix}
			\geq 0 
		\end{equation}
		for all $t \geq 0$.
\end{assumption}
The resulting dynamics can be described by 
\begin{equation} \label{eq:sys_dyn_G_hat_noisy}
    \begin{split}
        \Dot{\eta}(t)&=\hat{A}_G(\rho)\eta(t) + \hat{B}_G(\rho) u(t)+\hat{B}_G(\rho)e(t),\\
        y(t)&=\hat{C}_G(\rho)\eta(t), \\
        u(t)&=\nabla f (y(t)),
    \end{split}
\end{equation}
with initial conditions $\eta(0)=\eta_{0}$.
\begin{remark}
Note that although we split the term $$\hat{B}_G(\rho)u_n(t)=\hat{B}_G(\rho)u(t) + \hat{B}_G(\rho)e(t)$$ in \eqref{eq:sys_dyn_G_hat_noisy}, the exact value $u(t)=\nabla f (y(t))$ is not available to the agents and the terms are separated only for the analysis.
\end{remark}
The assumptions made so far ensure the existence of the equilibrium $\eta_*,\,u_*=0,\,e_*=0,\,y_*$ such that for any trajectory $\rho:[0,\infty) \rightarrow \mathcal{P}^N$,
\begin{equation} \label{eq:sys_dyn_G_hat_noisy_equilibrium}
    \begin{split}
        0&=\hat{A}_G(\rho)\eta_*, \\
        y_*&=\hat{C}_G(\rho)\eta_*, \\
        0&=\nabla f (y_*).
    \end{split}
\end{equation}
\begin{Problem} \label{prob:formation_control}
Under Assumptions \ref{assum:psi_Sml}, \ref{assm:path_to_leaders}, \ref{assm:equilibrium} and \ref{assum:noise}, derive sufficient conditions independent of $N$ under which the trajectories $(\eta(t),u(t),e(t),y(t),\rho(t))$ satisfying dynamics \eqref{eq:sys_dyn_G_hat_noisy} remain bounded and $y(t)$ converges exponentially (with a specified rate $\alpha$) to the minimizer $y_*$ of $f$.
\end{Problem}
\section{Review of Supporting Results}\label{sec:reviewSupportingResults}
Results from \cite{Datar.12102021,Datar.2022} are now reviewed in this section which support the main results presented in Section \ref{sec:main_results}.
Concretely, the main result presented in Section \ref{sec:main_results} assume that $f\in \mathcal{S}(m,L)$ and prove convergence to the minimizer of $f$ under this assumption.
We therefore review earlier results that characterize the minimizers of $f$ and give conditions to verify if $f$ is in $\mathcal{S}(m,L)$.
\subsection{Properties of $f$}\label{sec:prop_f}
We define two grounded Laplacians as in \cite{Datar.2022} that help in verifying $f\in  \mathcal{S}(m,L)$.
\begin{definition}{(\cite{Datar.2022})}\label{defn:Lgs_defn}
For a given graph $\mathcal{G}$ of order $N$ (with its corresponding Laplacian $\mathcal{L}$), a set of informed agents $\mathcal{V}_l$ and constants $0<m_{\psi}\leq L_{\psi}$, let the grounded Laplacians $\mathcal{L}_s$ and $\mathcal{L}_b$ be defined by
\begin{equation} \label{eq:Lgs_defn}
    \begin{split}
        \mathcal{L}_s=\mathcal{L}+m_{\psi} E \textnormal{ and }       \mathcal{L}_b=\mathcal{L}+L_{\psi} E, \\
    \end{split}
\end{equation}
where $E$ is a diagonal matrix with the $i^{th}$ diagonal entry equal to $1$ if $i \in \mathcal{V}_l$ and equal to $0$ otherwise.
\end{definition}
\begin{lemma}{(\cite{Datar.2022})} \label{lemm:f_in_Sml}
For a given graph $\mathcal{G}$ of order $N$ (with its corresponding Laplacian $\mathcal{L}$), a set of informed agents $\mathcal{V}_l$, a scalar field $\psi$, a formation reference vector $r$ and constants $m_{\psi}$, $L_{\psi}$ such that $0<m_{\psi}\leq L_{\psi}$, let $f$ be as defined in Definition \ref{defn:f} and $\mathcal{L}_s$, $\mathcal{L}_b$ be as defined in Definition \ref{defn:Lgs_defn}.
Then, for constants $m$, $L$ such that $0<m \leq L$, the following statements are equivalent:
\begin{enumerate}
    \item[1)] $f\in \mathcal{S}(m,L) \textnormal{ for all } \psi \in \mathcal{S}(m_{\psi},L_{\psi})$,
    \item[2)] $m I \preceq  \mathcal{L}_s$ and $ \mathcal{L}_b \preceq L I$.
\end{enumerate}
\end{lemma}
\begin{lemma}{(\cite{Datar.2022})}\label{lemm:m_in_Sml_exists}
Let $0<m_{\psi}\leq L_{\psi}$ be fixed.
There exist constants $m$, $L$ with $0<m\leq L$, such that $f \in \mathcal{S}(m,L)$ for every $\psi \in \mathcal{S}(m_{\psi},L_{\psi})$, if and only if the graph $\mathcal{G}$ and the set of informed agents $\mathcal{V}_l$ satisfy Assumption \ref{assm:path_to_leaders}.
\end{lemma}
\begin{lemma}{(\cite{Datar.2022})}\label{lemm:minimizers}
Let $\psi$ satisfy Assumption \ref{assum:psi_Sml}, $\mathcal{G}$ and $\mathcal{V}_l$ satisfy Assumption \ref{assm:path_to_leaders} and $f$ be as defined in Definition \ref{defn:f}.
Then the following statements hold:
\begin{enumerate}
    \item[1)] If $r=0$, then $z$ is the minimizer of $f$ if and only if $z=\mathbf{1}_N \otimes y_{\textnormal{opt}}$.
    \item[2)] If $\mathcal{V}_l=\{i\}$ for some $i \in \mathcal{V}$, then $z$ is the minimizer of $f$, if and only if $z_j=y_{\textnormal{opt}}+(r_j-r_i)$ for all $j \in \mathcal{V}$.
\end{enumerate}
\end{lemma}

\begin{remark}
Note that if $r$ is chosen such that $r_i=0$ for the informed agent $i \in \mathcal{V}_l$, then $r$ just encodes the desired positions of the agents with the coordinate system such that the source and the informed agent are located at the origin. 
\end{remark}
\begin{remark}\label{rem:CompObj}
Scenarios involving multiple informed agents and a non-zero formation reference $\hat{r}$ are difficult to characterize because the terms $\frac{1}{2}(y-r)^T\mathcal{L}_{(d)}(y-r)$ and $\sum_{i \in \mathcal{V}_l} \psi(y_i)$ have competing objectives.
The equilibrium can thus result in a situation where none of the informed agents are at the source and agents are not at desired relative distances.
\end{remark}
\begin{lemma}{(\cite{Datar.2022})}\label{lemm:minimizers_quadratic_radially_symmetric}
Let $\psi$ satisfy Assumption \ref{assum:psi_Sml}, $\mathcal{G}$ and $\mathcal{V}_l$ satisfy Assumption \ref{assm:path_to_leaders} and $f$ be as defined in Definition \ref{defn:f} with $z$ minimizing $f$.
Then $\{y: \nabla \psi(z_i)^T(z_i-y)\geq 0 \textnormal{ for all } i \in \mathcal{V}_l \}$ contains $y_{\textnormal{opt}}$.
Furthermore, the following stronger conclusions are obtained if $\psi$ is quadratic or radially symmetric.
\begin{enumerate}
    \item[1)] If $\psi$ is radially symmetric around the source, i.e., it has the form $\psi(y)=\psi_r(||y-y_{\textnormal{opt}}||)$ for some function $\psi_r : \mathbb{R} \rightarrow \mathbb{R}$, then the minimizer $z$ of $f$ is such that $y_{\textnormal{opt}}$ lies in the convex hull of $\{z_i : i \in \mathcal{V}_l\}$. 
    \item[2)] If $\psi$ is quadratic, i.e., it has the form $\psi(y)=y^TQy+b^Ty+c$, then, the minimizer $z$ of $f$ satisfies $y_{\textnormal{opt}}=\frac{1}{|\mathcal{V}_l|} \sum_{i \in \mathcal{V}_l} z_i$, i.e., the center of mass of informed agents is at the minimizer $y_{\textnormal{opt}}$ of $\psi$.
\end{enumerate}
\end{lemma}
\begin{remark}
Note that the conclusions in Lemma \ref{lemm:minimizers_quadratic_radially_symmetric} are obtained irrespectively of the choice of desired formation vector $r$. The key reason for this is that since we consider undirected graphs, we have equal and opposite interaction forces which cancel each other when deriving the center of mass dynamics. See \cite{Datar.2022} for details.
\end{remark}
\subsection{Zames-Falb $\alpha-$IQCs parameterization for $ \mathcal{S}(m,L)$}
Another central ingredient in the analysis presented in Section \ref{sec:main_results} is the $\alpha$-IQC result for $\mathcal{S}(m,L)$ which is summarized next.
Let $\Pi=\left[
\begin{array}{c|c}
	A_{\Pi} & B_{\Pi} \\
	\hline
	C_{\Pi} & D_{\Pi}
\end{array}
\right]$ and the set $\mathbb{P}$
be as defined in Appendix \ref{app:ZF_params}.
\begin{theorem}{(\cite{Datar.12102021})}\label{theom:theorem_LMI_ZF}
Let $\Tilde{u},\Tilde{y} \in \mathcal{L}_{2e}[0,\infty)$ be related by $\Tilde{u}=\nabla f(\Tilde{y}+y_*)$, where $f \in \mathcal{S}(m,L)$ and $y_*$ minimizes $f$.
Then, the signal $\Tilde{z}$ as defined by 
\begin{equation}\label{eq:signal_defn_z}
    \Tilde{z}(t)=\int_0^t C_{\Pi} e^{A_{\Pi}(t-\tau)}B_{\Pi}
\begin{bmatrix}
\Tilde{y}(\tau)\\
\Tilde{u}(\tau)
\end{bmatrix}d\tau
+ 
D_{\Pi}
\begin{bmatrix}
\Tilde{y}(t)\\
\Tilde{u}(t)
\end{bmatrix}
\end{equation}
satisfies for any $\alpha \geq 0$, the $\alpha-$IQC
\begin{equation}\label{eq:zPz_pos}
    \int_0^T e^{2\alpha t}\Tilde{z}^T(t) (P \otimes I) \Tilde{z}(t) dt \geq 0 \quad \forall P \in \mathbb{P},\forall T \geq 0.
\end{equation}
\end{theorem}
The LTI system $\Pi$ will be referred to as the ZF-multiplier and the order $\nu$ of the LTI system $\pi$ (See Appendix) used in the construction of $\Pi$ is called the order of the ZF multiplier.
\section{Main Analysis Results}\label{sec:main_results}
Building on the results presented in the last section (mainly Theorem \ref{theom:theorem_LMI_ZF}), the main results are derived in this section. Theorem \ref{theom:theorem_perf_analysis_noisy} presents a sufficient condition in the form of an LMI that guarantees exponential convergence with rate $\alpha$. 
The obtained LMI is then reduced in size to obtain an equivalent LMI in Lemma \ref{lemm:decomposition}  such that it is independent of the network size $N$. Finally, Theorem \ref{theom:theorem_perf_analysis_noisy_qlpv} extends these results to quasi-LPV systems.
\begin{theorem}\label{theom:theorem_perf_analysis_noisy}
	Let the multiplicative noise $e(t)$ in (\ref{eq:noise_condition}) satisfy Assumption \ref{assum:noise}.
	If there exists $\mathcal{X}\succ 0$, $P \in \mathbb{P}$, $\alpha >0$ and $\lambda\geq 0$ such that for all $\bar{\rho} \in \mathcal{P}^N$,
	\begin{align}\label{eq:perf_LMI_noise}
			&\begin{bmatrix}
				\mathcal{A}(\bar{\rho})^T\mathcal{X}+\mathcal{X}\mathcal{A}(\bar{\rho})+2\alpha \mathcal{X}  & (*) \\
				\mathcal{B}(\bar{\rho})^T\mathcal{X}    & \mathbf{0}
			\end{bmatrix} \nonumber \\
			&\hspace{1cm}+			
			(*)
			\begin{bmatrix}
				P \otimes I_{Nd} &\\&\lambda (M\otimes I_{Nd})
			\end{bmatrix}
			\begin{bmatrix}
				\mathcal{C}_1(\bar{\rho}) & \mathcal{D}_1\\
				\mathcal{C}_2 &\mathcal{D}_2\\
			\end{bmatrix}
			 \preceq 
			0,
	\end{align}
	where 
	\begin{align*}
	   &\left[\begin{array}{c|c}
		\mathcal{A}(\bar{\rho})    &  \mathcal{B}(\bar{\rho})\\
		\hline
		\mathcal{C}_1(\bar{\rho})    &  \mathcal{D}_1\\
		\mathcal{C}_2    &  \mathcal{D}_2\\
	\end{array}\right] \nonumber =\left[ 
	\begin{array}{cc}
		\Pi &\\
		&I_{2Nd}\\
	\end{array}\right]
	\left[\begin{array}{c|ccc}
		\hat{A}_G(\bar{\rho}) & \hat{B}_G(\bar{\rho}) & \hat{B}_G(\bar{\rho})\\		
		\hline
		\hat{C}_G(\bar{\rho})  & \mathbf{0}& \mathbf{0}\\
		\mathbf{0}  & I_{Nd}       & \mathbf{0}\\
		\hline
		\mathbf{0}  & \mathbf{0}& I_{Nd}\\
		\mathbf{0}  & I_{Nd}       & \mathbf{0}\\
	\end{array}
	\right]
	\end{align*}
	and $$M=\begin{bmatrix}-1&0\\0&\delta^2 \end{bmatrix},$$
	then, any trajectory $y$ generated from dynamics \eqref{eq:sys_dyn_G_hat_noisy} with model matrices satisfying Assumption \ref{assm:equilibrium}, $\rho:[0,\infty) \rightarrow \mathcal{P}^N$ and $f \in \mathcal{S}(m,L)$ converges to the minimizer $y_*$ of $f$ with rate $\alpha$, i.e., $\exists \kappa \geq 0$ such that  
	$||y(t)-y_*(t)||\leq \kappa e^{-\alpha t}$ holds for all $t\geq 0$.
\end{theorem}
\textit{Proof.}
Let $\xi=\begin{bmatrix}
	\eta - \eta_* \\
	x_{\pi}
\end{bmatrix},$ where $\eta$ is the state of $\hat{G}(\rho)$, $x_{\pi}$ is the state of the multiplier $\Pi$ and $\eta_*$ is the equilibrium state satisfying $\hat{A}_G\eta_*=\mathbf{0}$ and since $\nabla f (y_*)=0$, $\tilde{u}=u$ and therefore, $\tilde{e}=e$.
The dynamics of $\xi$ can be represented by
\begin{equation} \label{eq:sys_dyn_psi_GI_noise}
	\begin{split}
		\Dot{\xi}&=\mathcal{A}(\rho)\xi + \mathcal{B}(\rho)\begin{bmatrix}
			u\\e
		\end{bmatrix}, \quad\quad \xi(0)=[\Tilde{\eta}_0^T \quad \mathbf{0}]^T,\\
		\begin{bmatrix}
			\Tilde{z}\\ e\\u
		\end{bmatrix}&=\begin{bmatrix}
	\mathcal{C}_1(\rho)\\\mathcal{C}_2
	\end{bmatrix}\xi + \begin{bmatrix}
	\mathcal{D}_1\\\mathcal{D}_2
\end{bmatrix}\begin{bmatrix}
			u\\\tilde{e}
		\end{bmatrix},
	\end{split}
\end{equation}
where $\Tilde{z}$ can be obtained from \eqref{eq:signal_defn_z} with signals $\Tilde{u}$ and $\Tilde{y}$, where $\Tilde{y}$ is the output of $\hat{G}(\rho)$ for input $\Tilde{u}$. 
Furthermore, dynamics  \eqref{eq:sys_dyn_G_hat_noisy} imply that $\Tilde{u},\Tilde{y}$ satisfy $\Tilde{u}=\nabla f(\Tilde{y}+y_{\textnormal{opt}})$. 
Hence, Theorem \ref{theom:theorem_LMI_ZF} implies
\begin{equation} \label{eq:zPz_pos_noisy}
	\int_0^T e^{2\alpha t}\Tilde{z}^T(t) (P \otimes I_d) \Tilde{z}(t) dt \geq 0 \quad \forall P \in \mathbb{P},\forall T \geq 0.
\end{equation}
Define a storage function $V(\xi)=\xi^T \mathcal{X} \xi$.
Using \eqref{eq:perf_LMI_noise}, we get,
\begin{equation*}
	\begin{split}
		&\frac{d}{dt}(V(\xi(t)))+2\alpha V(\xi(t))\\
		&=(*)
		\begin{bmatrix}
			\mathcal{A}(\rho(t))^T\mathcal{X}+\mathcal{X}\mathcal{A}(\rho(t))+2\alpha \mathcal{X}  & \mathcal{X}\mathcal{B}(\rho(t))\\
			\mathcal{B}(\rho(t))^T\mathcal{X}    & \mathbf{0}
		\end{bmatrix}
		\begin{bmatrix}
			\xi(t)\\ \Tilde{u}(t)\\e(t)
		\end{bmatrix}\\ 
		&\leq		
		-(*)
		\begin{bmatrix}
				P \otimes I_{Nd} &\\&\lambda (M\otimes I_{Nd})
			\end{bmatrix}
		\begin{bmatrix}
			\mathcal{C}_1(\rho(t)) & \mathcal{D}_1\\
			\mathcal{C}_2&\mathcal{D}_2
		\end{bmatrix}
		\begin{bmatrix}
			\xi(t)\\ \Tilde{u}(t)\\e(t)
		\end{bmatrix} \\
		&=-\Tilde{z}^T(t) (P \otimes I_{Nd}) \Tilde{z}(t) - \lambda \begin{bmatrix}
			e^T&u^T
		\end{bmatrix}(M\otimes I_{Nd})  \begin{bmatrix}
		e\\u
	\end{bmatrix} \\
	&\leq-\Tilde{z}^T(t) (P \otimes I_{Nd}) \Tilde{z}(t),
	\end{split}
\end{equation*}
where the last inequality follows from Assumption \ref{assum:noise}.
Rearranging, multiplying by $e^{2\alpha t}$ and integrating from $0$ to $T$, we obtain,
\begin{equation*}\label{eq:diff_diss_ineq_noisy}
	\begin{split}
		\frac{d}{dt}(e^{2\alpha t}V(\xi(t))) + e^{2\alpha t}\Tilde{z}^T(t) (P \otimes I_{Nd}) \Tilde{z}(t) &\leq 0,\\
		 e^{2\alpha T}V(\xi(T)) + \int_0^T e^{2\alpha t}\Tilde{z}^T(t) (P \otimes I_{Nd}) \Tilde{z}(t) dt &\leq V(\xi(0)).
	\end{split}
\end{equation*}
Using \eqref{eq:zPz_pos} and $\mathcal{X}>0$, we get that $V(\xi(T)) \leq e^{-2\alpha T}V(\xi(0))$ implying $$||\Tilde{y}(T)|| \leq ||C_G|| \sqrt{\textnormal{cond}(\mathcal{X})} ||\xi(0)||  e^{-\alpha T}$$ for all $T\geq 0$.

\hfill $\square$

Owing to the structure of the multipliers used, the LMI \eqref{eq:perf_LMI_noise} is next decomposed to obtain an equivalent LMI independent of the network size $N$. 
\begin{lemma} \label{lemm:decomposition}
The following statements are equivalent:
\begin{enumerate}
    \item $\exists \mathcal{X}\succ 0,P \in \mathbb{P},\,\lambda\geq 0$ and $\alpha >0$ such that \eqref{eq:perf_LMI_noise} holds for all $\bar{\rho} \in \mathcal{P}^N$.
    \item $\exists \mathcal{X}_0\succ 0,P \in \mathbb{P},\,\lambda\geq 0$ and $\alpha >0$ such that for all $\tilde{\rho} \in \mathcal{P}$,
\begin{align}\label{eq:perf_LMI_noise_single}
			&\begin{bmatrix}
				\mathcal{A}_0(\tilde{\rho})^T\mathcal{X}_0+\mathcal{X}_0\mathcal{A}_0(\tilde{\rho})+2\alpha \mathcal{X}_0  & (*) \\
				\mathcal{B}_0(\tilde{\rho})^T\mathcal{X}_0    & \mathbf{0}
			\end{bmatrix} \nonumber \\
			&\hspace{0.1cm}+			
			(*)
			\begin{bmatrix}
				P \otimes I_d &\\&\lambda (M \otimes I_d)
			\end{bmatrix}
			\begin{bmatrix}
				\mathcal{C}_{10}(\tilde{\rho}) & \mathcal{D}_{10}\\
				\mathcal{C}_{20} &\mathcal{D}_{20}\\
			\end{bmatrix}
			 \preceq 
			0,
	\end{align}
		where 
	\begin{align*}
	   &\left[\begin{array}{c|c}
		\mathcal{A}_0(\tilde{\rho})    &  \mathcal{B}_0(\tilde{\rho})\\
		\hline
		\mathcal{C}_{10}(\tilde{\rho})    &  \mathcal{D}_{10}\\
		\mathcal{C}_{20}    &  \mathcal{D}_{20}\\
	\end{array}\right]= 
	\left[ 
	\begin{array}{cc}
		\Pi_0 & \\
		&I_{2d}\\
	\end{array}\right]
	\left[\begin{array}{c|ccc}
		A_G(\tilde{\rho}) & B_G(\tilde{\rho}) & B_G(\tilde{\rho})\\	
		\hline
		C_G(\tilde{\rho})  & \mathbf{0}& \mathbf{0}\\
		\mathbf{0}  & I_{d}       & \mathbf{0}\\
		\hline
		\mathbf{0}  & \mathbf{0}& I_{d}\\
		\mathbf{0}  & I_{d}       & \mathbf{0}\\
	\end{array}
	\right].
	\end{align*}
	\end{enumerate}
\end{lemma}
\textit{Proof}.
It can be shown that there exist permutation matrices $T_1$, $T_2$ such that
\begin{align*}
    \left[\begin{array}{c|c}
		\mathcal{A}(\bar{\rho})    &  \mathcal{B}(\bar{\rho})\\
		\hline
		\mathcal{C}_1(\bar{\rho})    &  \mathcal{D}_1\\
		\mathcal{C}_2    &  \mathcal{D}_2\\
	\end{array}\right]&=\left[\begin{array}{c|c}
		\hat{\mathcal{A}}_0(\bar{\rho})    &  \hat{\mathcal{B}}_0(\bar{\rho})\\
		\hline
		T_1\hat{\mathcal{C}}_{10}(\bar{\rho})    &  T_1\hat{\mathcal{D}}_{10},\\
		T_2\hat{\mathcal{C}}_{20}    &  T_2\hat{\mathcal{D}}_{20}\\
	\end{array}\right]\\
    T_1^T(P \otimes I_{Nd})T_1 & = I_N \otimes (P \otimes I_d), \\
    T_2^T(M \otimes I_{Nd})T_2&= (I_N \otimes (M \otimes I_d)).
\end{align*}
    
Substituting this in \eqref{eq:perf_LMI_noise}, an equivalent LMI is obtained as
	\begin{align}\label{eq:perf_LMI_noise_kroneckered}
			&\begin{bmatrix}
				\hat{\mathcal{A}}_0(\bar{\rho})^T\mathcal{X}+\mathcal{X}\hat{\mathcal{A}}_0(\bar{\rho})+2\alpha \mathcal{X}  & (*) \\
				\hat{\mathcal{B}}_0(\bar{\rho})^T\mathcal{X}    & \mathbf{0}
			\end{bmatrix} \nonumber \\
			&\hspace{0.1cm}+			
			(*)\begin{bmatrix}
				I_N \otimes(P \otimes I_d) &\\&I_N \otimes (\lambda M \otimes I_d)
			\end{bmatrix}
			\begin{bmatrix}
				\hat{\mathcal{C}}_{10}(\bar{\rho}) & \hat{\mathcal{D}}_{10}\\
				\hat{\mathcal{C}}_{20} &\hat{\mathcal{D}}_{20}\\
			\end{bmatrix}
			 \preceq 
			0.
	\end{align}
Now, let $\exists \mathcal{X}\succ 0,P \in \mathbb{P}$ such that \eqref{eq:perf_LMI_noise} holds for all $\rho \in \mathcal{P}^N$ which implies that \eqref{eq:perf_LMI_noise_kroneckered} holds.
Let $e_i$ denote the $i^{\textnormal{th}}$ canonical unit vector of dimension $2N$ and let $U_i=[e_i\,\,e_{N+i}] \otimes I_{n_x+2d}$.
Multiplying LMI \eqref{eq:perf_LMI_noise_kroneckered} from the left by $U_1^T$ and from the right by $U_1$, we get that LMI \eqref{eq:perf_LMI_noise_single} holds with $\mathcal{X}_0$ equal to the first principle minor $\mathcal{X}_{11}$ of $\mathcal{X}$.
Since $\mathcal{X} \succ 0$ implies $\mathcal{X}_{11} \succ 0$, it constitutes a feasible solution to LMI \eqref{eq:perf_LMI_noise_single} proving $(1) \implies (2)$.

Finally, let $\mathcal{X}=I_N \otimes \mathcal{X}_0$.
By defining the permutation matrix $U=[U_1\,U_2\,\cdots,U_N]$, LMI \eqref{eq:perf_LMI_noise_kroneckered} can be block-diagonalized by left and right multiplying  by $U^T$ and $U$, respectively, with each diagonal block equal to LMI \eqref{eq:perf_LMI_noise_single}.
Thus, if $\exists \mathcal{X}_0\succ 0,P \in \mathbb{P}$ such that for all $\tilde{\rho} \in \mathcal{P}$, \eqref{eq:perf_LMI_noise_single} is feasible, then $\mathcal{X}=I_N \otimes \mathcal{X}_0\succ 0$ solves \eqref{eq:perf_LMI_noise} proving $(2) \implies (1)$.

\hfill $\square$

\subsection{Extensions to q-LPV systems}
This section extends the results to quasi-LPV systems (see \cite{shamma1992linear}) of the form \eqref{eq:sys_dyn_G_hat_noisy} where the scheduling parameter trajectory is uniquely determined from the state of the system. 
Without going into the details, it is assumed that sufficiently smooth solutions exist for the considered quasi-LPV systems.
\begin{assumption} \label{assum:qlpv}
Assume there exists a function $g:\mathbb{R}^{n_{\eta}}\rightarrow\mathbb{R}^{n_{\rho}}$ such that $\rho(t)=g(\eta(t))$ holds for all $t\in[0,\infty)$ and with $\mathcal{P}^{-1}_g=\{\eta\in \mathbb{R}^{n_{\eta}}:g(\eta)\in \mathcal{P}^N\}$, there exists $c\geq 0$ such that 
\begin{equation*}
  \mathcal{B}(\eta_*,c)=\left\{ \eta:||\eta-\eta_*|| \leq c \right\} \subset \mathcal{P}^{-1}_g.
\end{equation*}
\end{assumption}

\begin{theorem}\label{theom:theorem_perf_analysis_noisy_qlpv}
Let the multiplicative noise $e$ satisfy Assumption \ref{assum:noise}.
Let there exist $\mathcal{X}_0\succ 0,P \in \mathbb{P},\,\lambda\geq 0$ and $\alpha >0$ such that for all $\tilde{\rho} \in \mathcal{P}$, LMI \eqref{eq:perf_LMI_noise_single} is satisfied.
Then for all initial conditions $\eta_0$ such that $||\eta_0-\eta_*||< \frac{c}{\sqrt{\textnormal{cond}(\mathcal{X})}}$, the trajectories generated from dynamics \eqref{eq:sys_dyn_G_hat_noisy} with model matrices satisfying Assumption \ref{assm:equilibrium}, $\rho(t)=g(\eta(t))$ and $f \in \mathcal{S}(m,L)$ converges to the minimizer $y_*$ of $f$ with rate $\alpha$, i.e., $\exists \kappa \geq 0$ such that  
	$||y(t)-y_*(t)||\leq \kappa e^{-\alpha t}$ holds for all $t\geq 0$.
\end{theorem}
\textit{Proof}.
From the proof of Lemma \ref{lemm:decomposition}, if $\mathcal{X}_0$ solves LMI \eqref{eq:perf_LMI_noise_single}, then $\mathcal{X}=I_N \otimes \mathcal{X}_0\succ 0$ solves \eqref{eq:perf_LMI_noise}.
Furthermore, if $\rho(t)=g(\eta(t)) \in \mathcal{P}^N$ for all $t \in [0,\infty)$, then Theorem \ref{theom:theorem_perf_analysis_noisy} completes the proof.
So, we just need to show that $\rho(t)=g(\eta(t)) \in \mathcal{P}^N$ for all $t \in [0,\infty)$, which, because of Assumption \ref{assum:qlpv} reduces to showing that $\eta(t) \in \mathcal{B}(\eta_*,c)$ for all $t \in [0,\infty)$. 
Note that since $\frac{c}{\sqrt{\textnormal{cond}(\mathcal{X})}}\leq c$, $\eta_0$ lies in the interior of $\mathcal{B}(\eta_*,c)$.
Let there exist some time $t_f> 0$ such that $||\eta(t_f)-\eta_*||=c$ and $||\eta(t)-\eta_*||< c$ for all $t\in[0,t_f)$.
Assumption \ref{assum:qlpv} then implies that $\rho(t)=g(\eta(t)) \in \mathcal{P}^N$ for all $t \in [0,t_f]$.
Theorem \ref{theom:theorem_perf_analysis_noisy} then implies that $||\eta(t_f)-\eta_*||^2\leq \textnormal{cond}(\mathcal{X}) ||\eta_0-\eta_*||^2<c^2=||\eta(t_f)-\eta_*||^2$, a contradiction.
Therefore, there exists no $t_f$ such that $||\eta(t_f)-\eta_*||=c$ and existence of continuous solutions then implies that $||\eta(t_f)-\eta_*||<c$ for all $t\in [0,\infty)$ completing the proof.
\hfill $\square$



\section{Numerical Results }\label{sec:numerical_results}
The code for producing the results and figures is available at \cite{IFAC2023code}.

\begin{remark}\label{rem:Bisection} Notice that the constraint \eqref{eq:perf_LMI_noise_single} in Lemma \ref{lemm:decomposition} is not linear in $\alpha$. A bisection search can be used to find the best $\alpha$ in \eqref{eq:perf_LMI_noise_single}, in which feasibility of \eqref{eq:perf_LMI_noise_single} is verified for different values of $\alpha$. This technique produces the curves in Fig. \ref{fig:quadrotor_optimal_kd_ZF_noisy}, \ref{fig:quadrotor_optimal_kd_ZF_noisy_conservatism}, \ref{fig:LvsAlpha_Delta_Kd20} and \ref{fig:KdvsAlpha_Delta_L70}. \end{remark}

\subsection{Single LTI Quadrotor}\label{eg:quadrotor}
\begin{figure}[h!]
	\centering
	\input{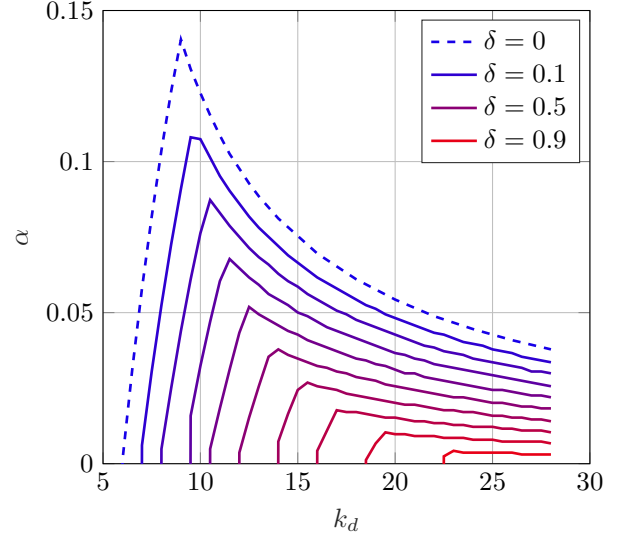}
	\caption{Performance estimates for quadrotor dynamics provided by first order ZF multipliers for $\psi \in \mathcal{S}(1,10)$ and varying $k_d$ with different bound $\delta \in \{0.1,0.2\cdots,0.9\}$ on the multiplicative noise.}
	\label{fig:quadrotor_optimal_kd_ZF_noisy}
\end{figure}

\begin{figure}[h!]
	\centering
%
%
\definecolor{mycolor1}{rgb}{0.09091,0.00000,0.90909}%
\scalebox{1}{
\begin{tikzpicture}

\begin{axis}[%
width=2.521in,
height=2.366in,
at={(0.758in,0.481in)},
scale only axis,
xmin=3,
xmax=30,
xlabel style={font=\color{white!15!black}},
xlabel={$k_d$},
ymin=0,
ymax=0.15,
ylabel style={font=\color{white!15!black}},
ylabel={$\alpha$},
yticklabels={0, 0, 0.05, 0.1, 0.15},
xmajorgrids,
ymajorgrids,
axis background/.style={fill=white},
title style={font=\bfseries},
title={},
legend style={legend cell align=left, align=left, draw=white!15!black}
]
\addplot [color=red!0!mycolor1, dashed, line width=1.0pt]
  table[row sep=crcr]{%
1	-1\\
1.5	-1\\
2	-1\\
2.5	-1\\
3	-1\\
3.5	-1\\
4	-1\\
4.5	-1\\
5	-1\\
5.5	-1\\
6	0.001220703125\\
6.5	0.0311279296875\\
7	0.0579833984375\\
7.5	0.081787109375\\
8	0.103759765625\\
8.5	0.123291015625\\
9	0.140380859375\\
9.5	0.130615234375\\
10	0.1226806640625\\
10.5	0.1153564453125\\
11	0.108642578125\\
11.5	0.1025390625\\
12	0.09765625\\
12.5	0.0927734375\\
13	0.0885009765625\\
13.5	0.0848388671875\\
14	0.0811767578125\\
16	0.069580078125\\
18	0.06103515625\\
20	0.0543212890625\\
22	0.048828125\\
24	0.0445556640625\\
26	0.0408935546875\\
28	0.037841796875\\
};
\addlegendentry{$\delta=0$}
\addplot [color=black, line width=1.0pt, only marks, mark=o, mark options={solid, black}]
  table[row sep=crcr]{%
5	-0.07095539474729\\
7	0.0581605673803028\\
9	0.140995625810792\\
11	0.10888760621175\\
13	0.0889845777293038\\
15	0.0753353299872735\\
17	0.0653601305878092\\
19	0.0577389631460918\\
21	0.051720900936853\\
23	0.046845506606932\\
25	0.0428140844217267\\
27	0.0394241088254656\\
29	0.0365332709136337\\
};
\addlegendentry{Ex. with $\delta=0$}

\addplot [color=red!50!mycolor1, line width=1.0pt]
  table[row sep=crcr]{%
1	-1\\
1.5	-1\\
2	-1\\
2.5	-1\\
3	-1\\
3.5	-1\\
4	-1\\
4.5	-1\\
5	-1\\
5.5	-1\\
6	-1\\
6.5	-1\\
7	-1\\
7.5	-1\\
8	-1\\
8.5	-1\\
9	-1\\
9.5	-1\\
10	-1\\
10.5	-1\\
11	-1\\
11.5	-1\\
12	0.003662109375\\
12.5	0.0152587890625\\
13	0.025634765625\\
13.5	0.035400390625\\
14	0.037841796875\\
14.5	0.0360107421875\\
15	0.0347900390625\\
15.5	0.0335693359375\\
16	0.0323486328125\\
16.5	0.03173828125\\
17	0.030517578125\\
17.5	0.0299072265625\\
18	0.0286865234375\\
18.5	0.028076171875\\
19	0.0274658203125\\
19.5	0.0262451171875\\
20	0.025634765625\\
20.5	0.0250244140625\\
21	0.0244140625\\
21.5	0.0238037109375\\
22	0.023193359375\\
22.5	0.0225830078125\\
23	0.02197265625\\
23.5	0.02197265625\\
24	0.0213623046875\\
24.5	0.020751953125\\
25	0.0201416015625\\
25.5	0.0201416015625\\
26	0.01953125\\
26.5	0.0189208984375\\
27	0.0189208984375\\
27.5	0.018310546875\\
28	0.018310546875\\
};
\addlegendentry{$\delta=0.5$}
\addplot [color=black, line width=1.25pt, draw=none, mark=asterisk, mark options={solid, black},only marks]
  table[row sep=crcr]{%
5	-0.199226365321045\\
7	-0.055972989506181\\
9	0.0401577268826836\\
11	0.0915509482192924\\
13	0.0732592079189861\\
15	0.0614669130668166\\
17	0.0548322584627636\\
19	0.0491797875201019\\
21	0.0442706506349314\\
23	0.0392183519500266\\
25	0.0363283142013289\\
27	0.0339632615215399\\
29	0.0308927840933636\\
};
\addlegendentry{Ex. with $\delta=0.5$}
\end{axis}

\end{tikzpicture}%
}
	\caption{Conservatism analysis for quadrotor dynamics provided by first order ZF multipliers for $\psi \in \mathcal{S}(1,10)$ and varying $k_d$ with different bound $\delta \in \{0.1,0.2\cdots,0.9\}$ on the multiplicative noise.}
	\label{fig:quadrotor_optimal_kd_ZF_noisy_conservatism}
\end{figure}
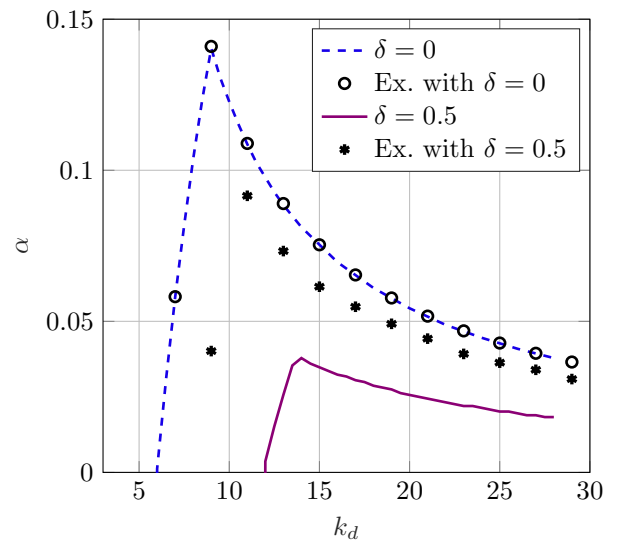

We consider an LTI example of a single linearized quadrotor model of order 12 with 6 states corresponding to positions and velocities and 6 states corresponding to orientation and angular velocities. 
We use a Linear-Quadratic-Regulator (LQR) based state-feedback controller tuned for zero steady-state error for tracking step position reference commands.
Let the noise level $\delta$ vary in the set $\{0.1,\cdots,0.9\}$ (see Assumption \ref{assum:noise}). 
The effect of varying $k_d$ (with fixed $k_p=1$) on the performance estimates for fixed allowable field set is shown in Fig. \ref{fig:quadrotor_optimal_kd_ZF_noisy}. 
It shows that the highest convergence rate estimate of  $0.14$ can be achieved for $k_d=9$ and demonstrates a method for tuning the gains for optimal convergence rate guarantees.
The trade-off between nominal performance and robustness against noise can be clearly seen.
The optimal values of $k_d$ increases as $\delta$ increases (thereby demanding more robustness) and the estimated performance therefore reduces.
Figure \ref{fig:quadrotor_optimal_kd_ZF_noisy_conservatism} investigates the conservatism by finding examples in the feasible set for which the convergence rate can be exactly computed.
Since we can find examples for the noise-free case that hit the estimates obtained from theory, there is no room for improvement in this case and there is no conservatism.
For the noisy case ($\delta=0.5$), however, there is a possible conservatism.
It is however emphasized that finding the worst case noise trajectory is not an easy task and analyzing conservatism in this case requires further efforts.

\subsection{Multiple Interacting Vehicles with Non-Linear Friction}\label{eg:unicycle}
Conider now a non-linear friction vehicle model 
\begin{equation}\label{eq:NLFsat}
m\ddot{\xpos}+b\vert\dot{\xpos}\vert\dot{\xpos} = u_F, 
\end{equation}
where $\xpos$ is the position, $m$ is the mass, $b$ is the viscous friction constant and $u_F$ is the force input.

A qLPV state space representation of \eqref{eq:NLFsat} is
\begin{equation}\label{eq:NLFssLPV}
    \begin{split}
        \begin{bmatrix} \dot{\xpos}\\ \dot{\xvel} \end{bmatrix} 
        &= 
        \begin{bmatrix} 0 & 1 \\ 0 & -\tfrac{b}{m}\rho \end{bmatrix} 
        \begin{bmatrix} \xpos\\ \xvel \end{bmatrix} 
        +
        \begin{bmatrix} 0 \\ \tfrac{1}{m} \end{bmatrix} 
        u_F,
        \qquad
        \\
        y
        &=
        \begin{bmatrix} 1 & 0 \end{bmatrix}
        \begin{bmatrix} \xpos\\ \xvel \end{bmatrix},
        \\
        \rho(t) 
        &=
        \vert \xvel(t)\vert,
    \end{split}
\end{equation}
where $\xvel$ is the velocity of the mass. 
A tracking output feedback controller is synthesized for a single vehicle using techniques presented in \cite{SCHERER2001361} to obtain a closed-loop tracking system.

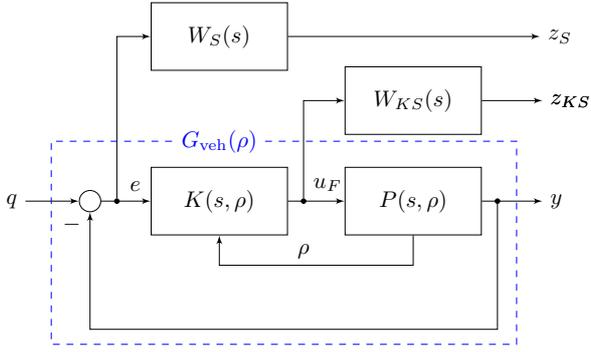
\begin{figure}[h!]
    \centering
     	\tikzstyle{block} = [draw, rectangle, 
    minimum height=3em, minimum width=6em]
\tikzstyle{sum} = [draw,  circle, node distance=1cm]
\tikzstyle{input} = [coordinate]
\tikzstyle{output} = [coordinate]
\tikzstyle{pinstyle} = [pin edge={to-,thin,black}]
\tikzstyle{joint} = [circle, draw, fill, inner sep=0pt, minimum size=2pt]

\scalebox{0.85}{
\begin{tikzpicture}[auto, node distance=2cm,>=latex']
    \node [input, name=input] {};
    \node [sum, right of=input] (sum) {};
    \node [block, right of=sum] (controller) {$K(s,\rho)$};
    \node [block, right of=controller,
            node distance=3cm] (system) {$P(s,\rho)$};
    \node [block, above=1.5cm of controller ] (WS) {$W_S(s)$};
    \node [block, above=0.5cm of system] (WKS) {$W_{KS}(s)$};
    
    \node [joint, right=2mm of sum] (joint_e) {};
    \node [joint, right=2mm of controller] (joint_u) {};
    \node [output, right of=system] (output) {};
    \node [output, right=3.9cm of WS] (outputS) {};
    \node [output, right of=WKS] (outputKS) {};
    \node [joint, right=2mm of system] (bifY){};

    \draw [->] (sum) --  (controller);
    \draw [->] (controller) --  (system);
    \draw [->] (input) -- (sum);
    \draw [->] (system) -- (output);
    \draw (system) -- ++(0,-1) [->] -| (controller);
    \draw (bifY) -- ++(0,-2) [->] -| (sum);
    \draw [->] (joint_e) |- (WS);
    \draw [->] (joint_u) |- (WKS);
    \draw [->] (WS) -- (outputS);
    \draw [->] (WKS) -- (outputKS);
 
    \node [above right=0mm and 0.5mm of joint_e] (e) {$e$};
    \node [above right=0mm and 0mm of joint_u] (u) {$u_F$};
    \node [right=0mm of outputKS] (zKS) {$z_{KS}$};
    \node [right=0mm of outputS] (zS) {$z_S$};
    \node [right=0mm of outputKS] (zKS) {$z_{KS}$};
    \node [left=0mm of input] (ry) {$q$};
    \node [right=0mm of output] (y) {$y$};
    \node [below=5mm of joint_u](LPVp){$\rho$};
    \node [below left=0mm and -1mm of sum] (minus) {$-$};
    \node[text=blue, above=0.1cm of controller] (Gveh) {$G_\text{veh}(\rho)$};
    \draw[blue, dashed] (Gveh.east)-|([xshift=-4mm]output.west)|-([yshift=-1.7cm]controller.south)-|([xshift=4mm]input.east)|-(Gveh.west);
        
\end{tikzpicture}
}
    \caption{Generalized LPV plant: closed loop $G_\text{veh}(\rho)$ with performance channel and shaping filters.}
    \label{fig:SKS_LPV}
\end{figure}

The LMIs \eqref{eq:perf_LMI_noise_single} are verified on the set $\mathcal{P}=[0,5]$ and it is emphasized that exponential performance is guaranteed only if the initial conditions are close enough to the source (see Theorem \ref{theom:theorem_perf_analysis_noisy_qlpv}). The overall system used for the design and synthesis of the LPV controller is 
shown in Fig. \ref{fig:SKS_LPV}, where a classic mixed sensitivity approach is
proposed to reduce the tracking error $e(t)$ and control effort $u_F(t)$, by shaping the $S(s)$ and $KS(s)$ sensitivity functions.

\begin{figure}[h!]
    \centering
	\input{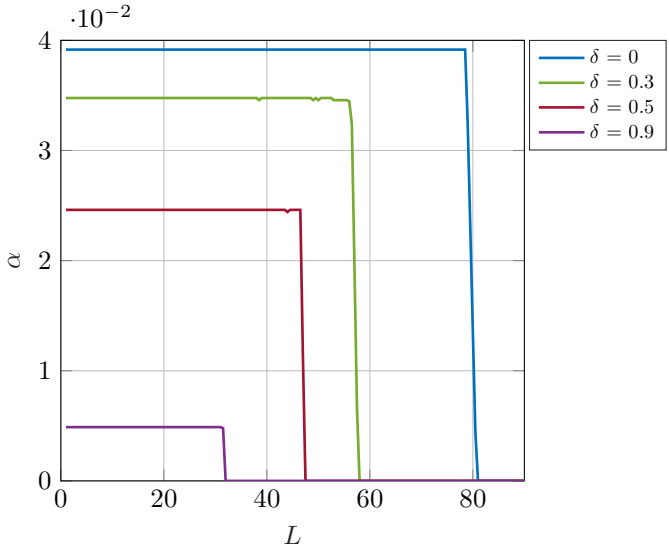}
    \caption{Variation of $L$ and $\delta$: convergence rate estimates for different $f \in \mathcal{S}(1,L)$ guaranteed by 5th ZF multipliers and $k_d=20$.}
    \label{fig:LvsAlpha_Delta_Kd20}
\end{figure}

\begin{figure}[h!]
    \centering
	\input{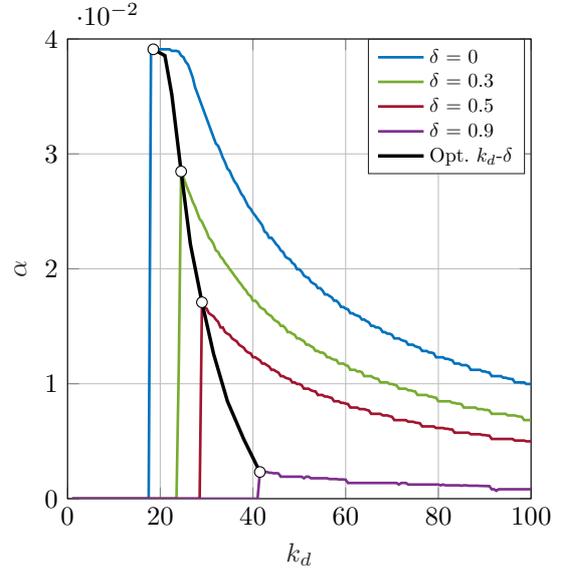}
    \caption{Variation of $k_d$ and $\delta$: convergence rates for different $k_d$ values and fixed $f \in \mathcal{S}(1,70)$ guaranteed by first order ZF multipliers.}
    \label{fig:KdvsAlpha_Delta_L70}
\end{figure}

The resulting synthesized controller is of 4th order, i.e. the overall closed loop has order 6 and $G_\text{veh}(\rho)$ from Fig. \ref{fig:SKS_LPV} is to be integrated into the framework of Fig. \ref{fig:plant G}. Numerical experiments in Fig. \ref{fig:LvsAlpha_Delta_Kd20} - \ref{fig:KdvsAlpha_Delta_L70} illustrate the influence of $L$, $k_d$ and $\delta$ on the convergence rate estimate $\alpha$ with 5th order ZF multipliers. 
Increasing the order of the multiplier reduces conservatism at the expense of computational cost. 
Lower order ZF multipliers produce conservative results and due to space restrictions are not included here.
Refer to \cite{Datar.2022} for a conservatism analysis of different types of multipliers.

Figure \ref{fig:LvsAlpha_Delta_Kd20} shows convergence rate estimates $\alpha$ for $f \in \mathcal{S}(1,L)$ and fixed $k_d=20$. Increasing $L$ enlarges the uncertainty set thereby giving a non-increasing estimate curve.
The intersection of the curves with the x-axis ($\alpha=0$) gives the maximum possible $L$ for which we can guarantee stability.
It can be seen that with noise levels of $0,\,0.3,\,0.5$ and $0.9$, stability is guaranteed up to $L=74,\,52,\,43$ and $27$ respectively, showing the reduction in robustness margins with increasing noise levels.
Furthermore, the estimated robust exponential performance rate $\alpha$ also reduces with increasing noise levels. 

Figures \ref{fig:LvsAlpha_Delta_Kd20} - \ref{fig:KdvsAlpha_Delta_L70} complement each other. In fact a 3d surface of $\alpha$ can be obtained for different ranges of $L$ and $k_d$ but would not be visually as clear as having them separately.

Similar to the Quadrotor example, Fig. \ref{fig:KdvsAlpha_Delta_L70} demonstrates the application of the methods to obtain controller gains with highest performance guarantees.
The optimal guaranteed performance $\alpha$ is achieved when $k_d\approx 20$ in the absence of noise ($\delta=0$); this optimal value becomes larger as the noise becomes stronger.
In other words, the damping gain $k_d$ needs to be increased to imporve robustness to noise.
Notice, that stability cannot be guaranteed for low damping, i.e. low $k_d$ values. 
The black curve can be used to pick the optimal $k_d$ based on a priori knowledge of the expected noise level $\delta$; this is illustrated in Fig. \ref{fig:KdvsAlpha_Delta_L70}. 

\begin{figure}[h!]
    \centering
    \input{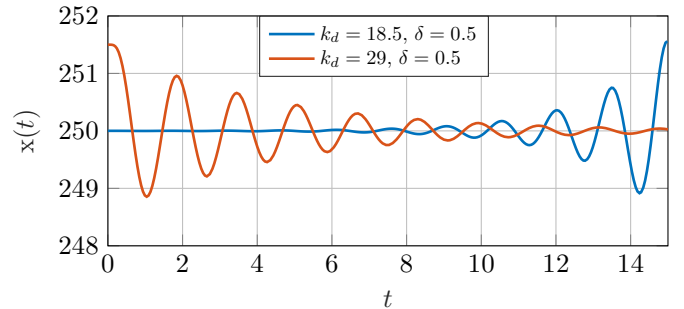}\\
    \caption{Single agent with noisy gradient and different values of $k_d$ optimized for $\delta=0$ and $\delta=0.5$ respectively.}
    \label{fig:1AgentNoise}
\end{figure}

Figure 
\ref{fig:1AgentNoise} depicts stable and unstable scenarios for a single agent attempting to locate the source of field $\psi \in \mathcal{S}(1,70)$ with $y_{\text{opt}}=250$ and noise level $\delta=0.5$. In blue, $k_d$ has been chosen to be optimal for noise-free gradients, i.e. $k_d=18.5$, which shows an unstable behaviour even if starting at $y(0)=250.001$; whereas the red curve converges to $y_\text{opt}$, since $k_d$ is optimized for noise levels of $\delta=0.5$ (see Fig. \ref{fig:KdvsAlpha_Delta_L70}).
The knowledge of the expected noise level can thus be incorporated into the design.

\begin{figure}[h!]
    \centering
    \input{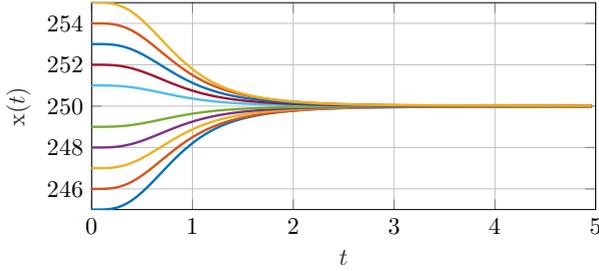}\\
    \caption{MAS scenario where 10 agents start at random initial conditions and reach the source of the field $\psi\in\mathcal{S}(1,66)$ with gradient noise $\delta= 0.5$ and $k_d=100$.}
    \label{fig:MAS}
\end{figure}

Finally, consider a scenario with multiple agents.
The external field $\psi$, the interconnection graph $\mathcal{G}$ and the set of informed agents $\mathcal{V}_l$ are restricted to be such that $f$ as defined in Definition \ref{defn:f} is in $\mathcal{S}(1,70)$.
The noise level is set to $\delta=0.5$.
Figure \ref{fig:KdvsAlpha_Delta_L70} can then be used to conclude exponential stability for $k_d\geq 29$. 
An example consensus scenario ($r_i=0$) with $10$ agents (all of them informed) interacting with a cycle graph, $\psi=\frac{1}{2}66(y-250)^2$ and $kd=100$ is considered which ensures that $f \in \mathcal{S}(1,70)$.
Figure \ref{fig:MAS} shows that agents reach consensus at the source.




\section{Conclusions and future work}\label{sec:conclusions}

A systematic approach for analyzing robust exponential performance of non-linear vehicles with qLPV representations is presented with a focus on cooperative source-seeking.
A priori information regarding the noise level can be explicitly considered in the analysis to trade-off exponential performance and robustness against noise, uncertain interconnections and scalar field. 
Numerical experiments demonstrate the potential use cases of this framework.
A more general analysis would include additive noise. Moreover, controller
synthesis within the presented framework is highly desirable.

\bibliography{rootIFAC}             

\appendix
\section*{APPENDIX}
\subsection*{Parameterization of ZF Multipliers}\label{app:ZF_params}
With $\beta=-1$, let
\begin{equation*}
	A_{\nu}=
	\begin{bmatrix} 
		\beta     & 0         &\dots     & 0 \\
		1           &\beta    & \ddots   &0  \\
		0           & \ddots    & \ddots   &0  \\
		\vdots      & 0         & 1        &\beta  
	\end{bmatrix},
	B_{\nu}=
	\begin{bmatrix} 
		1 \\
		0\\
		\vdots\\
		0 
	\end{bmatrix}.
\end{equation*}
Let $A_{\nu}^{\alpha}=A_{\nu}-2\alpha I$ and $\pi=\left[\begin{array}{c|c}
    A_{\nu}^{\alpha} & B_{\nu} \\
    \hline 
    0 & 1\\
    I_{\nu} & \mathbf{0}
\end{array}\right]$.
Finally, define 
\begin{align*}
\pi_{m,L}&=\begin{bmatrix}
    \pi&0\\0&\pi
\end{bmatrix}\begin{bmatrix}
    -m&1\\L&-1
\end{bmatrix}
=\left[
\begin{array}{cc|cc}
A_{\nu}^{\alpha}         &\mathbf{0}                         &-mB_{\nu}       &B_{\nu} \\
\mathbf{0}                          &A_{\nu}^{\alpha}        &LB_{\nu}        &-B_{\nu} \\
\hline
\mathbf{0}                          &\mathbf{0}                         &-m                      &1\\
I_{\nu}                  &\mathbf{0}                         &\mathbf{0}                 &\mathbf{0}\\
\mathbf{0}                          &\mathbf{0}                         &L                       &-1\\
\mathbf{0}                          &I_{\nu}                 &\mathbf{0}                 &\mathbf{0}
\end{array}
\right],
\end{align*}
$\Pi=\pi_{m,L} \otimes I_{Nd}$ and $\Pi_0=\pi_{m,L} \otimes I_{d}$.
With 
\begin{align*}
R_{\nu}&=\textnormal{diag}(\frac{1}{\sqrt{0!}},\cdots,\frac{1}{\sqrt{(\nu-1)!}}),\\
\left[
\begin{array}{c|c}
\Tilde{A}_{\nu} & \Tilde{B}_{\nu} \\
\hline
\Tilde{C}_{\nu} & \Tilde{D}_{\nu}
\end{array}
\right]&=\left[ 1 \quad \frac{s}{(s-\beta)^{\nu -1}} \quad \hdots \quad \frac{s^{\nu -1}}{(s-\beta)^{\nu -1}}\right]^T,
\end{align*}
consider the LMIs in variables $H \in \mathbb{R}$, $P_1, P_3 \in \mathbb{R}^{\nu}$, $\mathcal{X}_1, \mathcal{X}_3 \in \mathbb{S}^{\nu-1}$ 
\renewcommand{\theequation}{A.\arabic{equation}}
\begin{align}
    H+(P_1+P_3)A_{\nu}^{-1}B_{\nu} &\geq 0, \label{eq:L1_norm_constraint}
   \\
    (*)
    \begin{bmatrix}
    \mathbf{0}    & \mathcal{X}_i & \mathbf{0} \\
    \mathcal{X}_i & \mathbf{0}    & \mathbf{0} \\
    \mathbf{0} & \mathbf{0} & \textnormal{diag}(P_i)
    \end{bmatrix}
    \begin{bmatrix}
    I                       & \mathbf{0} \\
    \Tilde{A}_{\nu}         & \Tilde{B}_{\nu} \\
    R_{\nu}\Tilde{C}_{\nu}  & R_{\nu}\Tilde{D}_{\nu} \\
    \end{bmatrix}
    &\succ 0 \label{eq:positivity}.
\end{align}
Finally, define the set
\begin{equation*}
\mathbb{P}=\left\{
\begin{bmatrix}
	\mathbf{0}  & \mathbf{0} & 	H  & -P_3 \\
	\mathbf{0}  & \mathbf{0} &	-P_1^T    & \mathbf{0}\\
		*&*    & \mathbf{0}  & \mathbf{0} \\
		 *&*    & \mathbf{0}  & \mathbf{0} 
\end{bmatrix}:H, P_1, P_3\textnormal{ satisfy \eqref{eq:L1_norm_constraint},\eqref{eq:positivity}}\right\}.
\end{equation*}
\end{document}